\documentclass[12pt]{article}
\usepackage{amssymb, amsmath}
\usepackage{color}     

\newtheorem{thm}{Theorem}[section]

\newtheorem{cor}{Corollary}[section]

\renewcommand{\a}{\alpha}

\renewcommand{\o}{\omega}

\newcommand{\bb}{\begin{equation}}
\newcommand{\ee}{\end{equation}}
\newcommand{\bq}{\begin{eqnarray}}
\newcommand{\eq}{\end{eqnarray}}
\newcommand{\bqn}{\begin{eqnarray*}}
\newcommand{\eqn}{\end{eqnarray*}}%28 March 2014

\begin{document}
\title{  Notes on the Liouville  type problem for the stationary Navier-Stokes equations in $\Bbb R^3$  }
 \author{Dongho Chae\\
\ \\
Department of Mathematics\\
Chung-Ang University\\
 Seoul 156-756, Republic of Korea\\
email: dchae@cau.ac.kr}
\date{}
\maketitle
\begin{abstract}
In this paper we study the Liouville type problem for the stationary Navier-Stokes equations in $\Bbb R^3$.
We deduce an asymptotic formula for an integral involving the head pressure, $Q=\frac12 |v|^2 +p$,  and its derivative over  domains enclosed by level surfaces of $Q$.
This formula provides us  with                                                                                                                                                                                                                                                                                                                                                                                                                                                                                                                                                                                                                                                                                                                                                                                                                                                                                                                                                                                                                                                                                                                                                                                                                                                                                                                                                                                                                                                                                                                                                                                                                                                                                                                                                                                                                                                                                                                                                                                                                                                                                                                                                                                                                                                                                                                                                                                                                                                                                                                                                                                                                                                                                                                                                                                                                                                                                                                                                                                                                                                                                                                                                                                                                                                                                                                                                                                                                                                                                                                                                                                                                                                                                                                                                                                                                                                                                                                                                                                                                                                                                                                                                                                                                                                                                                                                                                                                                                                                                                                      new sufficient condition for the triviality of solution to the Navier-Stokes equations. 
  \noindent{\bf
keywords:} steady Navier-Stokes equations, Liouville type theorem \end{abstract}

\section{Introduction}
 \setcounter{equation}{0}

In this paper we are concerned on the Liouville type problem for the n-dimensional steady Navier-Stokes equations in $\Bbb R^3$.
\bb\label{sns}
\left\{ \aligned v\cdot \nabla v &=-\nabla p +\Delta v,\\
\nabla \cdot v&=0
\endaligned
\right.
\ee
equipped with the uniform decay condition at spatial infinity,
\bb
\label{bc1}
v(x)\to 0 \quad \text{as} \quad |x|\to 0 .
\ee
Here, $v=(v_1 (x), v_2 (x),  v_3 (x))$ is a vector field in $\Bbb R^3$, and $p=p(x)$ is a scalar field.  
Obviously  $(v,p)$ with $v=0$ and $p=$constant  is a trivial solution to \eqref{sns}-\eqref{bc1}. An important question is if there is other nontrivial solution.
This uniqueness problem, or  equivalently Liouville type problem is now hot issue in the community of mathematical fluid mechanics.
In general we impose an extra condition, the finiteness of the Dirichlet integral,
\bb\label{dir}
\int_{\Bbb R^3} |\nabla v|^2 dx <+\infty.
\ee
The Liouville type problem for the solution of \eqref{sns}-\eqref{bc1} together with \eqref{dir} is posed explicitly  in Galdi's book\cite[Remark X.9.4,pp.729]{gal}, where
the triviality of solution is deduced under stronger assumption $v\in L^{\frac92} (\Bbb R^3)$.
 % In the case of the n-dimensional Navier-Stokes system with $n\geq 4$  the problem for  \eqref{sns}-\eqref{dir}  is  resolved in \cite{gal} by similar argument to the case  of $v\in L^{\frac92} (\Bbb R^3)$, and the case $n=2$ is done by Gilbarg-Weinberger in \cite{gil}.  
%The case of $n=3$ is wide open, and 
There are other numerous partial results(see e.g. \cite{cha1,cha2, cha3, ser1, ser2, ser3, koz, cham} and the references therein) proving the triviality of solution  to \eqref{sns}-\eqref{dir} under various sufficient conditions. 
In the following theorem we derive an asymptotic blow-up rate for an integral of the head pressure over a region enclosed by level surfaces of the head pressure.
The integral has the same scaling property as the Dirichlet integral.
For $k\in \Bbb N$ let us define
$$ \exp_k (x)= \overbrace{\exp(\exp(\cdots (\exp x)\cdots  )}^{k-\text{times}}, \qquad  \log_k (x)= \overbrace{\log(\log(\cdots (\log x)\cdots  )}^{k-\text{times}},
$$
and set $ \exp_k (1):= e_k$, $e_0(x)=\log_0 (x):=1$ for all $x\in \Bbb R$. \begin{thm}
Let $(v,p)$ be a smooth solution of \eqref{sns} satisfying \eqref{bc1}-\eqref{dir}, and $Q=\frac12 |v|^2 +p$.  Then,
for all  $k\in \Bbb N \cup \{ 0\}$ we have
\bb\label{th1}
\lim_{\lambda \to 0} (\log_{k+1}(1/\lambda))^{-1} \int_{\{ |Q|>\lambda \}} \frac{ |\nabla Q|^2 }{|Q| \prod_{j=0}^k \log_j \left(\frac{ e_k |m|}{|Q|}\right)} dx=  \int_{\Bbb R^3} |\o|^2 dx\ee
Therefore, if there exists $k\in \Bbb N\cup\{0\}$ such that 
\bb\label{th2}
 \int_{\{ |Q|>\lambda \}} \frac{ |\nabla Q|^2 }{|Q| \prod_{j=0}^k \log_j \left(\frac{ e_k |m|}{|Q|}\right)} dx=  o\left(   \log_{k+1}(1/\lambda) \right)
\ee
as $\lambda \to 0$, then $v=0$ on $\Bbb R^3$.
\end{thm}
In order to appreciate the above theorem we consider the following simple case with  $k=0$. 
\begin{cor}
Under the the same assumptions  as in Theorem 1.1 there holds 
\bb\label{co1}
\lim_{\lambda \to 0} (\log(1/\lambda))^{-1} \int_{\{ |Q|>\lambda \}} \frac{ |\nabla Q|^2 }{|Q| } dx=  \int_{\Bbb R^3} |\o|^2 dx, \ee
and  if 
\bb\label{co2}
\int_{\{ |Q|>\lambda \}} \frac{ |\nabla Q|^2 }{|Q| } dx=  o\left(   \log(1/\lambda) \right)
\ee
as $\lambda \to 0$, then $v=0$ on $\Bbb R^3$. 
\end{cor}
\noindent{\bf Remark 1.1 } We note  that $\int_{\Bbb R^3} |\nabla \sqrt{|Q|} |^2 dx = \frac14 \int_{\Bbb R^3 } \frac{ |\nabla Q|^2 }{|Q| } dx$ has the same scaling as $\int_{\Bbb R^3} |\nabla v|^2 dx$.
As an immediate consequence of the above corollary we have the triviality, $v=0$ if $\int_{\Bbb R^3 } \frac{ |\nabla Q|^2 }{|Q| } dx<+\infty$. Indeed, if $Q=$constant, then \eqref{ns}  below implies $\o=0$, and this together with $\nabla \cdot v=0$
leads to  $v=\nabla h$ for a harmonic function $h$, which is zero from the condition \eqref{bc1} by the Liouville theorem for a harmonic function. It would be interesting to note the following inequality, which follows from the Sobolev embedding,
\bb \label{re1}
\left(\int_{\Bbb R^3} |Q|^3 dx\right)^{\frac16}  \leq C\left( \int_{\Bbb R^3} |\nabla \sqrt{|Q|} |^2 dx\right)^{\frac12}.
\ee
We see that $v$ is trivial if the right hand side  of  \eqref{re1} is finite. The left hand side  of \eqref{re1} is, however,  easily shown to be finite from the hypothesis of Theorem 1.1. Indeed, using the well-known Calderon-Zygmund inequality\cite{ste} for the pressure, $\Delta p=-\sum_{j,k=1}^3\partial_j\partial_k (v_jv_k)$,
we have
\begin{align}\label{re2}
&\left(\int_{\Bbb R^3} |Q|^3 dx\right)^{\frac16} \leq C\left(\int_{\Bbb R^3}|v|^6  dx\right)^{\frac16} +  C\left(\int_{\Bbb R^3}|p|^3  dx\right)^{\frac16} \cr
&\qquad \qquad\le C\left(\int_{\Bbb R^3} |v|^6 dx\right) ^{\frac16} \le C \left(\int_{\Bbb R^3} |\nabla v|^2 dx \right)^{\frac12} <+\infty.
\end{align}
\section{Proof of Theorem 1.1}
 \setcounter{equation}{0}
%\noindent{\bf Proof } 
Given smooth solution $(v,p)$ of \eqref{sns} we  set the head pressure $Q = p +\frac12 |v |^2$ and the vorticity tensor $\o=\nabla \times v $. Then,  it is well-known that the following holds
\bb\label{ns}
\Delta Q -v \cdot \nabla  Q =|\o|^2.
\ee
Indeed,  multiplying the first equation of \eqref{sns} by $v$, we have
\bb\label{nns1}
v\cdot \nabla Q=\Delta \frac{|v|^2}{2} -|\nabla v|^2.
\ee
Taking divergence of \eqref{sns}, we are led to
\bb\label{nns2}
0=\Delta p+ \rm{Tr}(\nabla v (\nabla v)^\top) .
\ee
Adding \eqref{nns2} to \eqref{nns1}, and observing $ |\nabla v|^2-   \rm{Tr}(\nabla v (\nabla v)^\top)=|\o|^2$,  we obtain \eqref{ns}.
We also recall that it is known (see e.g. Theorem X.5.1, p. 688 [2]) that  the condition \eqref{bc1} together with \eqref{dir} implies that  $ p(x) \to \bar{p}$  for some constant $ \bar{p}$ 
 as $|x| \to +\infty$.
Therefore, re-defining $Q -\bar{p}$ as  new head pressure, we may assume 
\bb\label{bc2}
Q(x)\to 0 \qquad \text{as}\quad |x|\to+\infty.
\ee
In view of  \eqref{bc2},  applying the maximum principle  to  \eqref{ns}, we find $Q\leq 0$ on $\Bbb  R^3$.
Moreover, by the maximum principle again, either $Q(x) = 0$  for all $x\in \Bbb R^3$, or $Q(x) < 0$ for all 
$x\in \Bbb R^3$.  Indeed, any point  $x_0 \in \Bbb R^3$ such that $Q(x_0) =0$ is a point of local maximum, which is not allowed unless $Q \equiv 0$ by the maximum principle.  Let $Q(x) \not=0$ on $\Bbb R^3$, then without the loss of generality we may assume $Q(x)<0$ for all $x\in \Bbb R^3$. 
We set 
$$\inf_{x\in \Bbb R^3} Q(x) =m<0. $$
Let us consider a function $f: \Bbb R_+ \to \Bbb R$, which is  Lipshitz continuous, and satisfying $f(0)=0.$  Given $\lambda \in (0, |m|)$, 
we multiply  \eqref{ns} by $f(-Q) $, and then integrating it over  the set $  \{x\in \Bbb R^3\, |\, |Q(x)|>\lambda \}$, we have
\begin{align}
\label{nss}
\int_{ \{|Q|>\lambda \}}  |\o|^2 f(-Q)  dx&= \int_{ \{|Q|>\lambda \}} \Delta Q  f(-Q)  dx\cr
 &\qquad\qquad-\int_{ \{|Q|>\lambda \}}  v\cdot \nabla Q f(-Q) dx \cr
 &:=I+J .
\end{align}
Let us define $F(s)=\int_0 ^s f(\sigma) d\sigma$.  Using the divergence theorem,  we compute
\begin{align}
\label{nsa}
J&= \int_{ \{|Q|>\lambda \}}   v\cdot \nabla  F(-Q) dx=\int_{ \{|Q|>\lambda \}}  \nabla \cdot \left\{v F(-Q)\right\} dx\cr
&= F(-\lambda ) \int_{ \{Q=-\lambda \}} v\cdot \nu dS = F(-\lambda ) \int_{\{ |Q|>\lambda \}} \nabla\cdot v dx=0,
\end{align}
where $\nu$ denotes the outward unit normal vector on the level surface.
  \begin{align}
  \label{nsb}
I&=\int_{ \{|Q|>\lambda \}}  \nabla \cdot \left\{\nabla Q f(-Q)\right\}dx+ \int_{ \{|Q|>\lambda \}}   |\nabla Q|^2   f'(-Q) dx  \cr
&:= I_1 +I_2.
\end{align}
By the divergence theorem and  \eqref{ns} we  have
\begin{align}
\label{nsc}
I_1&= f(-\lambda ) \int_{\{ Q=-\lambda \}} \frac{\partial Q}{\partial \nu} dS  = f(-\lambda )\int_{ \{|Q|>\lambda \}}  \Delta Q dx\cr
&= f(-\lambda)\int_{ \{|Q|>\lambda \}} |\o|^2 dx,
 \end{align}  
 where in the third equality we used \eqref{ns} and the fact
\begin{align}\label{12}
 &\int_{ \{|Q|>\lambda \}}  v\cdot \nabla Q dx= \int_{ \{|Q|>\lambda \}} \nabla \cdot (v Q) dx \cr
 &\qquad= -\lambda \int_{ \{ Q=-\lambda \} } v\cdot  \nu dS= -\lambda \int_{ \{|Q|>\lambda \}} \nabla \cdot v  dx =0.
\end{align}
Summarizing the above computation,  we have
\bb\label{11}
 \int_{ \{|Q|>\lambda \}}  |\o|^2( f(|Q|)  -f(\lambda )) dx=\int_{ \{ |Q|>\lambda\}}  |\nabla Q|^2   f'(|Q|) dx .
\ee 
Given $\alpha<0$, choosing $f(s)= (\log_k (\frac{e_k |m|}{s} ) )^{\alpha}$ in \eqref{11},  we have
 
   \begin{align}\label{12}
&\int_{ \{ |Q|>\lambda \}} |\o|^2\left\{\left(\log_k \frac{ e_k |m|}{\lambda} \right)^{\a} -\left(\log_k \frac{e_k |m|}{|Q|}\right)^{\a}\right\}dx\cr
&\qquad= \a \int_{\{ |Q|>\lambda \}} \frac{ |\nabla Q|^2 }{|Q| \prod_{j=0}^k \log_j \left(\frac{ e_k |m|}{|Q|}\right)} \left(\log_k\frac{e_k|m|}{ |Q|}\right)^{\a}  dx .
 \end{align}    

  Let us set $\a =-\beta (\log_{k+1} (1/\lambda))^{-1}$, $\beta >0$ in \eqref{12}, and pass $\lambda \to 0$. Then, we deduce
 \begin{align}\label{13}
 &\lim_{\lambda \to 0} (\log_{k+1} (1/\lambda ))^{-1} \int_{\{ |Q|>\lambda \}} \frac{ |\nabla Q|^2 }{|Q| \prod_{j=0}^k \log_j \left(\frac{ e_k |m|}{|Q|}\right)} \left(\log_k\frac{e_k|m|}{ |Q|}\right)^{-\beta (\log_{k+1} (1/\lambda))^{-1}}  dx\cr
& \qquad= \frac{1}{\beta} (1-e^{-\beta}) \int_{\Bbb R^3} |\o|^2 dx,
 \end{align}
 where we used  the dominated convergence theorem together with the fact that 
\begin{align}\label{14}
&\lim_{\lambda \to 0}\left(\log_k \frac{e_k |m|}{|Q|}\right)^{-\beta (\log_{k+1} (1/\lambda))^{-1}} =1, \cr
 & \lim_{\lambda \to 0}  \left(\log_k \frac{ e_k |m|}{\lambda} \right)^{-\beta (\log_{k+1} (1/\lambda))^{-1}} = e^{-\beta}.
\end{align}
Note that 
\bb\label{15}
\left(\log_k\frac{e_k|m|}{ \lambda }\right)^{-\beta (\log_{k+1} (1/\lambda))^{-1}} < \left(\log_k\frac{e_k|m|}{ |Q|}\right)^{-\beta (\log_{k+1} (1/\lambda))^{-1}} \le 1  
\ee
on $ \{ x\in \Bbb R^3\, |\, |Q(x)| >\lambda \}$.  Inserting \eqref{15} into \eqref{13}, and using \eqref{14}, we  obtain
  \begin{align*}
   & e^{-\beta} \lim_{\lambda \to 0} (\log_{k+1} (1/\lambda ))^{-1} \int_{\{ |Q|>\lambda \}} \frac{ |\nabla Q|^2 }{|Q| \prod_{j=0}^k \log_j \left(\frac{ e_k |m|}{|Q|}\right)}  dx  \cr
   &\qquad  \le \frac{1}{\beta}  (1- e^{-\beta}  ) \int_{  \Bbb R^3}  |\o|^2dx\cr
   &\qquad \le \lim_{\lambda \to 0} (\log_{k+1} (1/\lambda ))^{-1} \int_{\{ |Q|>\lambda \}} \frac{ |\nabla Q|^2 }{|Q| \prod_{j=0}^k \log_j \left(\frac{ e_k |m|}{|Q|}\right)}dx,  \cr   
       \end{align*}
  which implies that
      \begin{align}\label{16}
   &    \frac{1}{\beta} (1-e^{-\beta})  \int_{\Bbb R^3} |\o|^2 dx \le        \lim_{\lambda \to 0} (\log_{k+1} (1/\lambda ))^{-1} \int_{\{ |Q|>\lambda \}} \frac{ |\nabla Q|^2 }{|Q| \prod_{j=0}^k \log_j \left(\frac{ e_k |m|}{|Q|}\right)}dx \cr
   &\qquad  \le     \frac{1}{\beta} (e^\beta -1)  \int_{\Bbb R^3} |\o|^2 dx.
       \end{align}  
  Passing $\beta \to 0$ in \eqref{16}, we obtain \eqref{th1}. $\square$\\
  \ \\
  
        $$\mbox{\bf Acknowledgements}$$
The author was partially supported by NRF grants 2016R1A2B3011647.
 

\begin{thebibliography}{1}
\bibitem{cha1}D. Chae, {\it Liouville-type theorem for the forced Euler equations and the Navier-
Stokes equations.}  Commun. Math. Phys., {\bf 326}  pp. 37-48,  (2014).
\bibitem{cha2}  D. Chae and T. Yoneda, {\it On the Liouville theorem for the stationary Navier-Stokes
equations in a critical space,} J. Math. Anal. Appl. {\bf 405},  no. 2, 706-710,  (2013).
\bibitem{cha3} D. Chae and J. Wolf, {\it On Liouville type theorems for the steady Navier- Stokes
equations in $\Bbb R^3$}, J. Diff. Eqns., {\bf 261}, pp. 5541-5560,(2016).
\bibitem{cham}D. Chamorro, O. Jarrin and P- G. Lemari\'{e}-Rieusset, {\it Some Liouville theorems for stationary
Navier-Stokes equations in Lebesgue and
Morrey spaces,} arXiv preprint. no. 1806.03003.
%\bibitem{esc}L. Escauriaza, G. Seregin and V. \v{S}ver\'ek, {\it $L_{3, \infty}$-solutions of the Navier-stokes equations and backward uniqueness,} 
%Russian Math. Surveys, {\bf 58}, pp. 211-250, (2003),
\bibitem{gal} G. P. Galdi, {\it An introduction to the mathematical theory of the Navier- Stokes
equations: Steady-State Problems.}  Second edition. Springer Monographs in
Mathematics. Springer, New York,  xiv+1018 pp, (2011).
\bibitem{gil} D. Gilbarg and H.F. Weinberger, {\it 
Asymptotic properties of steady plane solutions of the Navier-Stokes equations with bounded Dirichlet integral, } Ann. Sc. Norm. Super. Pisa (4) 5,
pp. 381-404, (1978).
%\bibitem{koch} G. Koch, N.  Nadirashvili, G.  Seregin, V. \v{S}ver\'ek, {\it Liouville theorems for the
%Navier-Stokes equations and applications, } Acta Mathematica, {\bf 203}, pp. 83-105 (2009).
\bibitem{koz} H. Kozono, Y. Terasawa, Y. Wakasugi, {\it A remark on Liouville-type theorems
for the stationary Navier-Stokes equations in three space dimensions,}  J.  Func. Anal, {\bf 272},  pp.804-818, (2017).
\bibitem{ser1} G. Seregin, {\it Liouville type theorem for stationary Navier-Stokes equations,} Nonlinearity,
{\bf 29},   pp. 2191-2195, (2016),
\bibitem{ser2} G. Seregin, {\it Remarks on Liouville type theorems for steady-state Navier-Stokes
equations,}  Algebra i Analiz, {\bf 30},  no.2,  pp. 238-248, (2018).
\bibitem{ser3} G. Seregin and W. Wang, {\it Sufficient conditions on Liouville type theorems for the 3D steady Navier-Stokes equations,} arXiv preprint no. 1805.02227.
\bibitem{ste}E. M. Stein, {\it Singular Integrals and Differentiability Properties of Functions,} Princeton Univ. Press. (1970).


\end{thebibliography}
\end{document}